\documentclass[final,5p,twocolumn]{elsarticle}%
\usepackage{amssymb, amsmath,accents}%\usepackage{pdfsync} % Inversesearch
\journal{Systems \& Control Letters}
\usepackage{overpic}
\usepackage{float}

\usepackage{hyperref}  % make clickable links in LaTeX :-)

\usepackage{xcolor}
\usepackage{pdfsync}

\usepackage{graphicx}
\graphicspath{{figures/}}
\usepackage{epsfig}
\usepackage{amsmath,accents,dsfont}
\usepackage{amssymb}
\usepackage{amsthm}
\usepackage{mathtools}
\usepackage{amsfonts}
\usepackage{lineno}
\usepackage{tikz}
\usetikzlibrary{positioning,shapes}
\usetikzlibrary{arrows}

\newcommand{\esssup}{\operatorname*{ess\;sup}}%

\allowdisplaybreaks

\newtheorem{theorem}{Theorem}
\newtheorem{definition}[theorem]{Definition}

\newtheorem{lemma}[theorem]{Lemma}
\newtheorem{assumption}[theorem]{Assumption}
\newtheorem{proposition}[theorem]{Proposition}
\newtheorem{remark}[theorem]{Remark}

\newcommand{\tm}{\times}%

\usepackage{csquotes}  % For \enquotes

\def\field#1{\mathbb #1}%
\def\R{\field{R}}%
\newcommand{\Rn}[1][n]{\R^{#1}}

\newcommand{\ep}{\varepsilon}%
\newcommand{\N}{\mathbb{N}}%
\newcommand{\KC}{\mathcal{K}}%
\newcommand{\VC}{\mathcal{V}}%
\newcommand{\LC}{\mathcal{L}}%
\newcommand{\UC}{\mathcal{U}}%
\newcommand{\PC}{\mathcal{P}}%
\newcommand{\rmd}{\mathrm{d}}%
\newcommand{\rmD}{\mathrm{D}}%
\newcommand{\inner}{\mathrm{int}}%
\DeclareMathOperator{\id}{id}

\def\K{\mathcal{K}}%
\def\LL{\mathcal{L}}%
\def\KL{\mathcal{KL}}%
\def\Kinf{\mathcal{K}_\infty}%
\newcommand{\Rp}{\R_+}%
\let\ol=\overline%
\let\ul=\underline%
\def\ellinf{\ell_\infty}%
\def\ellinfp{\ell_\infty^+}%
\def\MHAF{\mathrm{MHAF}}%

\newcommand{\dist}{\mathrm{dist}}%

\newcommand{\unit}{\mathds{1}}%

\usepackage{xspace}	

 \title{\LARGE ISS small-gain criteria for infinite networks with linear gain functions}
\tnotetext[mytitlenote]{
The work of A.~Mironchenko is supported by the DFG through the grant MI 1886/2-1. The work of N. Noroozi is supported by the DFG through the grant WI 1458/16-1. 
The work of C. Kawan is supported by the German Research Foundation (DFG) through the grant ZA 873/4-1. The work of M. Zamani is supported in part by the DFG through the grant ZA 873/4-1 and the H2020 ERC Starting Grant AutoCPS (grant agreement No.~804639).
}
 
 \author[passauaddress]{Andrii~Mironchenko}
\address[passauaddress]{Faculty of Computer Science and Mathematics, University of Passau, 94032 Passau, Germany}

\author[lmuaddress]{Navid Noroozi\corref{mycorrespondingauthor}}
\cortext[mycorrespondingauthor]{Corresponding author (Email address: navid.noroozi@lmu.de)}

\author[lmuaddress]{Christoph Kawan}

\author[cubaddress,lmuaddress]{Majid Zamani}
\address[lmuaddress]{Institute of Informatics, LMU Munich,  Germany}
\address[cubaddress]{Computer Science Department, University of Colorado Boulder,USA}

\begin{document}

\begin{abstract}

\vspace{-0.3cm}
This paper provides a Lyapunov-based small-gain theorem for input-to-state stability (ISS) of networks composed of infinitely many finite-dimensional systems.
We model these networks on infinite-dimensional $\ell_{\infty}$-type spaces.
A crucial assumption in our results is that the internal Lyapunov gains, modeling the influence of the subsystems on each other, are linear functions. Moreover, the gain operator built from the internal gains is assumed to be subadditive and homogeneous, which covers both max-type and sum-type formulations for the ISS Lyapunov functions of the subsystems.
As a consequence, the small-gain condition can be formulated in terms of a generalized spectral radius of the gain operator.
Through an example, we show that the small-gain condition can easily be checked if the interconnection topology of the network has some kind of symmetry.
While our main result provides an ISS Lyapunov function in implication form for the overall network, an ISS Lyapunov function in a dissipative form is constructed under mild extra assumptions.

\end{abstract}


\maketitle

%\mir{We should decide between 'spectral small-gain condition' and 'spectral radius condition'. The first one is used only by us, and the second one by many researchers, but the meaning of the expression varies from the context.}

\section{Introduction}

%\mir{I do not understand 'poorly scale well' below.}

The paradigm shift towards a hyper-connected world composed of a large number of smart networked systems calls for developing scalable tools for rigorous analysis and synthesis of dynamical networks.
Standard tools in analysis and design of control systems, however, scale poorly with the size of such networks.
%Neglecting this fact will give rise to a system whose performance and stability indices dramatically degrade with the increase of the size of the system.
This calls for the development of new approaches ensuring the independence of stability/performance indices from the network size. A promising approach to address the scalability issue is to over-approximate a large-but-finite network with an infinite network. The latter refers to a network of infinitely many (finite-dimensional) subsystems and can be viewed as the limit case of the former in terms of the number of participating subsystems. Having studied an infinite network, it can be shown that the performance/stability indices achieved for the infinite network are transferable to any finite truncation of the infinite network~\cite{NMW.20,Barooah.2009,BPD02}.

Small-gain theory is widely used for analysis and design of finite-dimensional networks. In particular, the integration of Lyapunov functions with small-gain theory leads to Lyapunov-based small-gain theorems, where each subsystem satisfies a so-called input-to-state stability (ISS) Lyapunov condition~\cite{Liu.2014,Dashkovskiy.2010}. If the gain functions associated to the ISS Lyapunov conditions fulfill a small-gain condition, ISS of the overall network can be concluded.
%; see, e.g.,~\cite[Thm.~5.3]{Dashkovskiy.2010}.%

Advances in infinite-dimensional ISS theory, e.g.~\cite{MiP20, KaK19,ZhZ18}, created a firm basis for the development of Lyapunov-based small-gain theorems for infinite networks, recently studied in, e.g.,~\cite{KMS19,NMK21,DaP20,DMS19a}. Particularly, it has been observed that existing small-gain conditions for finite networks do not ensure ISS of an infinite network, even if all subsystems are linear~\cite{DMS19a}.%

To address such issues, we recently developed a Lyapunov-based small-gain theorem in a summation form (called sum-type small-gain theorem) in~\cite{KMS19}, i.e., the overall ISS Lyapunov function is a linear combination of ISS Lyapunov functions of the subsystems.
The small-gain condition, proposed in~\cite{KMS19}, is tight in the sense that it cannot be relaxed under the assumptions imposed on the network, cf.~\cite[Sec.~VI.A]{KMS19} for more details.
The overall state space (as well as the overall input space) is, however, modeled as an $\ell_p$-type space with $1 \leq p < \infty$.
This condition on the state space requires that each state vector has a vanishing tail, which clearly fails for several applications. To remove this condition, one needs to model the overall state space as an $\ell_{\infty}$-type space. 
%The sum-type small-gain condition in~\cite{KMS19} is not well-defined for an $\ell_{\infty}$-type space. \textcolor{red}{[Not so clear what this means...]}
This motivates the introduction of small-gain theorems in a maximum form (called max-type small-gain theorem).
In this type of small-gain theorem, each component of the gain operator, encoding the information on the interaction between subsystems, is associated with its corresponding subsystem and is expressed as a maximum over gain functions associated to its own neighboring subsystems.
In~\cite{DaP20,NMW.20}, max-type small-gain conditions are provided for continuous-time and discrete-time infinite networks, respectively.
Nevertheless, in both works, the gain functions are assumed to be uniformly less than the identity function. To let gain functions be larger than the identity, the so-called robust strong small-gain condition has been introduced in~\cite{DMS19a} and a method to construct a so-called \emph{path of strict decay} was accordingly proposed.
Building upon~\cite{Dashkovskiy.2010}, this path of strict decay is then used to construct an ISS Lyapunov function for the overall network.
However, in~\cite{DMS19a}, the existence of a linear path of strict decay is assumed, which is very restrictive.
This assumption is removed in~\cite{KMZ21}, where ISS of an infinite network was shown provided that there is a nonlinear path of strict decay with respect to the gain operator. Furthermore, in \cite{KMZ21} sufficient conditions for the existence of a nonlinear path of decay are provided.
In particular, it is shown that a (nonlinear) path of decay exists if the gain operator fulfills a robust strong small-gain condition and the discrete-time system induced by the gain operator is uniformly globally asymptotically stable.

All of the above small-gain results developed in an $\ell_{\infty}$-type space are given in a pure maximum formulation.
Clearly, a pure max formulation is not necessarily the best choice, in general, as one might need to add conservatism in the calculation of ISS Lyapunov bounds. For instance, one might unnecessarily upper-bound terms expressed in summation by those in maximum.
Inspired by~\cite{Dashkovskiy.2010}, our work introduces a new, more general small-gain theorem (cf.~Theorem~\ref{thm_mainres}), where the gain operator is monotone, subadditive and homogeneous of degree one. Our small-gain theorem covers both sum and max formulations of small-gain theorems as special cases. Therefore, our setting treats sum and max formulations in a unified, generalized way. To obtain such a formulation, we assume that the internal gain functions are linear. Our small-gain condition is expressed in terms of a spectral radius condition for the gain operator.%

%\mir{You say only about formulations, and new conditions. Please do not be shy and say, that we show a new small-gain theorem, and please refer to it by number.}

Our small-gain theorem relies on Proposition~\ref{prop:eISS-criterion-homogeneous-systems} below.
In this proposition, we show that the \emph{small-gain condition is equivalent to uniform global exponential stability of the monotone discrete-time system induced by the gain operator}, which is of great significance on its own, and extends several known criteria for linear discrete-time systems summarized in~\cite{GlM20}.
Moreover, \emph{this proposition provides a linear path of strict decay through which we construct an ISS Lyapunov function for the overall network.} 
For max and sum formulations, the spectral radius condition (i.e., the small-gain condition) admits explicit formulas for which a graph-theoretic description is provided.
Via an example of a spatially invariant network, the computational efficiency of these conditions is illustrated. To the best of our knowledge, our small-gain theorem and Proposition~\ref{prop:eISS-criterion-homogeneous-systems} are also novel for finite-dimensional systems and recover~\cite[Thm.~2]{Lu05} and~\cite[Lem.~1.1]{Rue10}, respectively, in max and sum formulations.%

%\mir{Next sentence is a copypaste from abstract.}
Though we mainly focus on the construction of an ISS Lyapunov function in an implication form for the overall system, in Proposition~\ref{prop:Nonlinear-scaling} we show that an ISS Lyapunov function in a dissipative form can be constructed under some mild extra assumptions.
In contrast to ODE systems \cite[Rem.~2.4]{SoW95}, it is unclear whether, in general, the existence of an ISS Lyapunov function in an implication form implies the existence of one in a dissipative form~\cite[Open Problem 2.17]{MiP20}.
Our result addresses this problem for Lyapunov functions constructed via small-gain design.

\vspace{-0.3cm}
\section{Technical setup} \label{sec:Technical setup}
\vspace{-0.2cm}
\subsection{Notation}

We write $\N$, $\R$ and $\Rp$ for the sets of positive integers, real numbers and nonnegative real numbers, respectively. We also write $\N_0 := \N\cup \{0\}$. Elements of $\Rn$ are by default regarded as column vectors. If $f:\R^n \rightarrow \R$ is differentiable at $x\in\R^n$, $\nabla f(x)$ denotes its gradient which is regarded as a row vector. By $\ellinf$ we denote the Banach space of all bounded real sequences $x = (x_i)_{i\in\N}$ with the norm $\|x\|_{\ellinf} := \sup_{i\in\N}|x_i|$. The subset $\ellinfp := \{ (x_i)_{i\in \N} \in \ellinf : x_i \geq 0,\ \forall i \in \N \}$ is a closed cone with nonempty interior, and is called the \emph{positive cone} in $\ellinf$. It induces a partial order on $\ellinf$ by ``$x \leq y$ if and only if $y - x \in \ellinfp$'' (which simply reduces to $x_i \leq y_i$ for all $i\in\N$). In any metric space $(X,d)$, we write $\inner(A)$ for the interior of a set $A$, $B_{\delta}(x)$ for the open ball of radius $\delta$ centered at $x$, and $\dist(x,A) := \inf_{y\in A}d(x,y)$ for the distance of a point $x$ to a set $A$. By $C^0(X,Y)$ we denote the set of all continuous functions from $X$ to (another space) $Y$. A function $f:\R_+ \rightarrow X$ is called \emph{piecewise right-continuous} if there are pairwise disjoint intervals $I_1 = [0,a_1)$, $I_2 = [a_1,a_2)$, $I_3 = [a_2,a_3),\ldots$ whose union equals $\R_+$, such that $f_{|I_n}$ is continuous for each $n \in \N$.%

A continuous function $\alpha:\Rp \to \Rp$ is positive definite, denoted by $\alpha\in\PC$, if $\alpha(0) = 0$ and $\alpha(r) > 0$ for all $r>0$. For the sets of comparison functions $\K$, $\Kinf$, $\LL$ and $\KL$, we refer to~\cite{Kellett.2014}. A function $\mu:\ellinfp \to [0,\infty]$ is called a \emph{monotone, homogeneous aggregation function} if it has the following properties:%
\begin{enumerate}
%\item\label{itm:positivity} Positivity: $\mu(s)>0$ for all $s\in %\inner(\ellinfp)$. \textcolor{red}{[Is positivity ever used?]}%
\item Homogeneity of degree one: $\mu(cs) = c\mu(s)$ for all $s\in \ellinfp$ and $c \geq 0$ (with the convention that $0 \cdot \infty = 0$).%
\item Monotonicity: $\mu(r) \leq \mu(s)$ for all $r,s\in \ellinfp$ s.t. $r \leq s$.%
\item Subadditivity: $\mu(r + s) \leq \mu(r) + \mu(s)$ for all $r,s\in \ellinfp$ (with the convention that $\infty + \infty = \infty$).%
\end{enumerate}
The set of all monotone, homogeneous aggregation functions is denoted by $\MHAF$. The concept of monotone aggregation functions has been introduced in~\cite{Dashkovskiy.2010}, and here it is adapted to our setting. Typical examples are $\mu(s) = \sum_{i\in\N} \alpha_i s_i$, $\mu(s) = \sup_{i\in\N}\alpha_i s_i$, or $\mu(s) = \max_{1\leq i \leq N} \alpha_i s_i + \sum_{i=N+1}^\infty \alpha_i s_i$, for some $(\alpha_i)_{i\in\N} \in\ellinfp$ with $\alpha_i\geq0$ for all $i\in\N$, and $N\in\N$.%

\subsection{Interconnections}

Consider a family of control systems of the form%
\begin{equation}\label{eq_subsystem}
  \Sigma_i:\quad \dot{x}_i = f_i(x_i,\bar{x}_i,u_i),\quad i \in \N.%
\end{equation}
This family comes with sequences $(n_i)_{i\in\N}$ and $(m_i)_{i\in\N}$ of positive integers as well as finite (possibly empty) sets $I_i \subset \N$, $i \notin I_i$, such that the following assumptions are satisfied:%
\begin{itemize}
\setlength\itemsep{0.1em}
\item The \emph{state vector} $x_i$ is an element of $\R^{n_i}$.%
\item The \emph{internal input vector} $\bar{x}_i$ is composed of the state vectors $x_j$, $j \in I_i$, and thus is an element of $\R^{N_i}$, where $N_i := \sum_{j \in I_i}n_j$.%
\item The \emph{external input vector} $u_i$ is an element of $\R^{m_i}$.%
\item The \emph{right-hand side} $f_i$ is a continuous function $f_i:\R^{n_i} \tm \R^{N_i} \tm \R^{m_i} \rightarrow \R^{n_i}$.%
\item For every initial state $x_{i0} \in \R^{n_i}$ and all essentially bounded inputs $\bar{x}_i(\cdot)$ and $u_i(\cdot)$, there is a unique (local) solution of the Cauchy problem $\dot{x}_i(t) = f_i(x_i(t),\bar{x}_i(t),u_i(t))$, $x_i(0) = x_{i0}$, which we denote by $\phi_i(t,x_{i0},\bar{x}_i,u_i)$.%
\end{itemize}
For each $i\in\N$, we fix (arbitrary) norms on the spaces $\R^{n_i}$ and $\R^{m_i}$, respectively. For brevity, we avoid adding indices to these norms, indicating to which space they belong, and simply write $|\cdot|$ for each of them. The interconnection of systems $\Sigma_i$, $i\in\N$, is defined on the state space%
\begin{equation*}
  X := \ell_{\infty}(\N,(n_i)) \!:=\! \{ x = (x_i)_{i\in\N} : x_i \in \R^{n_i},\ \sup_{i \in \N}|x_i| < \infty \}.%
\end{equation*}
This space is a Banach space with the $\ell_{\infty}$-type norm $\|x\|_X := \sup_{i\in\N}|x_i|$. The space of admissible external input values is likewise defined as the Banach space $U := \ell_{\infty}(\N,(m_i))$, $\|u\|_U := \sup_{i\in\N}|u_i|$. Finally, the class of admissible external input functions is 
\begin{equation*}
\UC := \{ u \in L_{\infty}(\R_+,U) : u \mbox{ is piecewise right-continuous} \},
\end{equation*}
which is equipped with the $L_{\infty}$-norm $\|u\|_{\UC} := \esssup_{t\in\R_+} |u(t)|_U$. Finally, the right-hand side of the interconnected system is defined by%
\begin{equation*}
  f:X \tm U \rightarrow \prod_{i\in\N}\R^{n_i},\ f(x,u) := (f_i(x_i,\bar{x}_i,u_i))_{i\in\N}.%
\end{equation*}
Hence, the interconnected system can formally be written as the following differential equation:%
\begin{equation}\label{eq_Sigma}
  \Sigma:\quad \dot{x} = f(x,u).%
\end{equation}
To make sense of this equation, we introduce an appropriate notion of solution. For fixed $x^0 \in X$ and $u \in \UC$, a function $\lambda:J \rightarrow X$, where $J$ is an interval of the form $[0,T)$ with $0 < T \leq \infty$, is called a \emph{solution} of the Cauchy problem $\dot{x}(t) = f(x(t),u(t))$, $x(0) = x^0 \in X$ provided that $s \mapsto f(\lambda(s),u(s))$ is a locally integrable $X$-valued function (in the Bochner integral sense) and $\lambda(t) = x_0 + \int_0^t f(\lambda(s),u(s))\, \rmd s$ for all $t \in J$. Sufficient conditions for the existence and uniqueness of solutions can be found in~\cite[Thm.~III.2]{KMS19}. We say that the system $\Sigma$ is \emph{well-posed} if local solutions exist and are unique. In this case, it holds that for all $t \in J$%
\begin{equation*}
  \pi_i(\phi(t,x^0,u)) = \phi_i(t,x^0_i,\bar{x}_i,u_i),%
\end{equation*}
where $\pi_i:X \rightarrow \R^{n_i}$ denotes the canonical projection onto the $i$-th component, $\bar{x}_i(\cdot) = (\pi_j(\phi(\cdot,x,u)))_{j\in I_i}$, and $x_i^0,u_i$ denote the $i$-th components of $x^0$ and $u$, respectively, see \cite[Sec.~III]{KMS19}.%

Throughout this paper, we assume that $\Sigma$ is well-posed and also has the so-called \emph{boundedness-implies-continuation (BIC) property} \cite{KaJ11b}. The BIC property requires that any bounded solution, defined on a compact time interval, can be extended to a solution defined on a larger time interval. If $\Sigma$ is well-posed, we call $\Sigma$ \emph{forward-complete} if every solution can be extended to $\R_+$.%

\subsection{Input-to-state stability}

We study input-to-state stability of $\Sigma$ via the construction of ISS Lyapunov functions.
% Definitions of these notions are given below.%

\begin{definition}\label{def:iss}
A well-posed system $\Sigma$ is called \emph{input-to-state stable (ISS)} if it is forward complete and there are $\beta\in\KC\LC$ and $\kappa\in\KC$ such that for all $x\in X$ and $u \in \UC$,%
\begin{equation*}
  \|\phi(t,x,u)\|_X \leq \beta(\|x\|_X,t) + \kappa(\|u\|_{\UC}) \mbox{\quad for all\ } t \geq 0.%
\end{equation*}
\end{definition}

A sufficient condition for ISS is the existence of an ISS Lyapunov function~\cite[Thm.~1]{DaM13}.%

\begin{definition}\label{def:implication-form}
A function $V:X \rightarrow \R_+$ is called an \emph{ISS Lyapunov function (in an implication form)} for $\Sigma$ if it satisfies the following conditions:%
\begin{enumerate}[(i)]
\item $V$ is continuous.%
\item There exist $\psi_1,\psi_2 \in \KC_{\infty}$ such that%
\begin{equation}\label{eq_isslf_coercivity}
  \psi_1(\|x\|_X) \leq V(x) \leq \psi_2(\|x\|_X) \mbox{\quad for all\ } x \in X.%
\end{equation}
\item There exist $\gamma\in\KC$ and $\alpha\in\PC$ such that for all $x\in X$ and $u\in\UC$ the following implication holds:%
\begin{equation}\label{eq_Lyap_impl}
  V(x) > \gamma(\|u\|_{\UC}) \quad \Rightarrow \quad \rmD^+ V_u(x) \leq -\alpha(V(x)),%
\end{equation}
where $\rmD^+ V_u(x)$ denotes the right upper Dini orbital derivative, defined as%
\begin{equation*}
  \rmD^+ V_u(x) := \limsup_{t \rightarrow 0^+} \frac{V(\phi(t,x,u)) - V(x)}{t}.%
\end{equation*}
\end{enumerate}
The functions $\psi_1,\psi_2$ are also called \emph{coercivity bounds}, $\gamma$ is called a \emph{Lyapunov gain}, and $\alpha$ is called a \emph{decay rate}.%
\end{definition}

To construct an ISS Lyapunov function $V$ for $\Sigma$, we follow a bottom-up approach in the sense that we exploit the interconnection structure and build $V$ from ISS Lyapunov functions of subsystems $\Sigma_i$ under an appropriate small-gain condition. We thus make the following assumption on each subsystem $\Sigma_i$.%

\begin{assumption}\label{ass_subsystem_iss}
For each subsystem $\Sigma_i$, there exists a continuous function $V_i:\R^{n_i} \rightarrow \R_+$, which is $C^1$ outside of $x_i=0$, and satisfies the following properties:%
\begin{enumerate}
\item[(L1)] There exist $\psi_{i1},\psi_{i2} \in \Kinf$ such that%
\begin{equation}\label{eq_subsystem_iss_coerc}
  \psi_{i1}(|x_i|) \leq V_i(x_i) \leq \psi_{i2}(|x_i|) \mbox{\quad for all\ } x_i \in \R^{n_i}.%
\end{equation}
\item[(L2)] There exist $\mu_i\in\MHAF$, $\gamma_{ij} \in \R_+$ ($j\in\N$), where $\gamma_{ij} = 0$ for all $j \notin I_i$, and $\gamma_{iu} \in \KC$, $\alpha_i \in \PC$ such that for all $x = (x_j)_{j\in\N} \in X$ and $u = (u_j)_{j\in\N} \in U$ the following implication holds:%
\begin{align}\label{eq_subsystem_orbitalder_est}
\begin{split}
  V_i(x_i) &> \max\Bigl\{\mu_i \big((\gamma_{ij} V_j (x_j)\big)_{j\in\N}),\gamma_{iu}(|u_i|) \Bigr\} \\
	&\Rightarrow \nabla V_i(x_i) f_i(x_i,\bar{x}_i,u_i) \leq -\alpha_i(V_i(x_i)) .
\end{split}
\end{align}
\end{enumerate}
The numbers $\gamma_{ij}$ (which are identified with the linear functions $r \mapsto \gamma_{ij}r$) are called \emph{internal gains}, while the functions $\gamma_{iu}$ are called \emph{external gains}.%
\end{assumption}

Given MHAFs $\mu = (\mu_i)_{i\in\N}$ and $\gamma_{ij}$ as in~\eqref{eq_subsystem_orbitalder_est}, we introduce the gain operator $\Gamma_\mu:\ell_{\infty}^+ \rightarrow \ell_{\infty}^+$ by%
\begin{equation}\label{eq:Gamma}
  \Gamma_\mu(s) := \Big(\mu_i \big((\gamma_{ij} s_j)_{j\in\N}\big) \Big)_{i\in\N}   \quad\forall s = (s_i)_{i\in\N} \in \ell_{\infty}^+.%
\end{equation}

This operator covers linear, max-linear and mixed types of gain operators. With the choice of all $\mu_i$ in~\eqref{eq:Gamma} either as summation or supremum, we have the following two special cases of the operator $\Gamma_\mu$:%
\begin{align}
  \Gamma_\oplus (s) &:= \Big(\sum_{j\in\N} \gamma_{ij}s_j \Big)_{i\in\N} \mbox{\quad for all\ } s = (s_i)_{i\in\N} \in \ell_{\infty}^+, \label{eq:Gamma-sum} \\
  \Gamma_\otimes	 (s) &:= \Big(\sup_{j\in\N} \gamma_{ij}s_j \Big)_{i\in\N} \mbox{\quad for all\ } s = (s_i)_{i\in\N} \in \ell_{\infty}^+. \label{eq:Gamma-max}  
\end{align} 

%\begin{assumption}\label{ass_gainop_wd}
%There is $\gamma_{\max} > 0$ such that%
%\begin{equation*}
%  \gamma_{ij} \leq \gamma_{\max} \mbox{\quad for all\ } i,j \in \N.%
%\end{equation*}
%\end{assumption}

The following assumption guarantees that $\Gamma_\mu$ is well-defined and induces a discrete-time system with bounded finite-time reachability sets, cf.~\cite{MiW18b}.%

\begin{assumption}\label{ass_gainop_wd}
The operator $\Gamma_\mu$ is well-defined, which is equivalent to
\begin{equation}\label{eq:uniformity}
  \|\Gamma_\mu(\unit)\|_{\ell_{\infty}} < \infty.%
\end{equation}
\end{assumption}

\begin{remark}\label{rem:well-posedness}
The condition \eqref{eq:uniformity} has a simple expression in terms of the gains $\gamma_{ij}$ if the $\mu_i$ have a particular structure. More specifically, if $\Gamma = \Gamma_\oplus$ is linear, then condition~\eqref{eq:uniformity} becomes $\sup_{i\in\N}\sum_{j\in\N} \gamma_{ij} < \infty$. If $\Gamma=\Gamma_\otimes$, condition~\eqref{eq:uniformity} reduces to $\sup_{i,j\in\N} \gamma_{ij} < \infty$.
\end{remark}

\begin{definition}\label{def:homogeneity-subadditivity-etc} 
An operator $\Gamma:\ell_\infty^+\to\ell_\infty^+$ is called
\begin{enumerate}
\item\label{itm:eISS-criterion-homog-0} monotone, if $\Gamma(r) \leq \Gamma(s)$  for all $r,s\in \ell_\infty^+$ s.t.  $r \leq s$.%
\item\label{itm:eISS-criterion-homog-1} homogeneous of degree one, if $\Gamma(cs) = c\Gamma(s)$ for all $s\in \ellinfp$ and $c \geq 0$.%
\item\label{itm:eISS-criterion-homog-2} subadditive, if $\Gamma(r + s) \leq \Gamma(r) + \Gamma(s)$ for all $r,s\in \ell_\infty^+$.%
\end{enumerate}
\end{definition}

The properties of the monotone aggregation functions $\mu_i$ immediately imply the statement of the next lemma.

\begin{lemma}\label{lem:Properties-Gamma-mu} 
The operator $\Gamma_\mu$ is monotone, subadditive and homogeneous of degree one.
\end{lemma}

Finally, we make the following uniform boundedness assumption on the external gains, which we need to prove our small-gain theorem.%

\begin{assumption}\label{ass_gammaiu_max}
There is $\gamma_{\max}^u \in \KC$ such that $\gamma_{iu} \leq \gamma_{\max}^u$ for all $i \in \N$.%
\end{assumption}

We develop a small-gain theorem ensuring ISS of the overall network $\Sigma$ from the assumptions imposed on subsystems $\Sigma_i$.
The small-gain theorem is given in terms of a spectral radius condition.
In Section~\ref{sec:equivalences}, we establish several characterizations of this condition. Then, in Section~\ref{sec:Small-gain theorem}, we present our small-gain theorem. While our setting mainly considers ISS Lyapunov functions in an implication form, in Section~\ref{sec:Dissipative formulation}, we show that the overall Lyapunov function $V$ is, modulo scaling, also an ISS Lyapunov function in a dissipative form under mild extra assumptions.

\section{Characterizations of the stability of the gain operator}\label{sec:equivalences}

Consider an operator $\Gamma:\ell_\infty^+\to\ell_\infty^+$ and the corresponding induced system%
\begin{equation}\label{eq_monotone_system-SGT}
  x(k+1) = \Gamma(x(k)),\quad k \in \N_0.%
\end{equation}

Based on~\cite[Prop.~16, 18]{MKG20}, we have the following characterizations for uniform exponential stability of the system~\eqref{eq_monotone_system-SGT} induced by a monotone, subadditive and homogeneous of degree one operator.
These characterizations play a crucial role in the verification of our small-gain condition, used in Theorem~\ref{thm_mainres} below.%

\begin{proposition}\label{prop:eISS-criterion-homogeneous-systems}
Let $\Gamma:\ell_{\infty}^+ \rightarrow \ell_{\infty}^+$ be an operator which is continuous, monotone, subadditive, homogeneous of degree one, and satisfies Assumption~\ref{ass_gainop_wd} (with $\Gamma$ instead of $\Gamma_\mu$). Then the following statements are equivalent:%
\begin{enumerate}[(i)]
\item The spectral radius condition%
\begin{equation}\label{eq_def_spectralradius}
  r(\Gamma) := \lim_{n \rightarrow \infty} \sup_{s \in \ell_\infty^+\atop \|s\|_{\ell_{\infty}} = 1}\|\Gamma^n(s)\|_{\ell_{\infty}}^{1/n} < 1.%
\end{equation}
\item \emph{Uniform global exponential stability (UGES)} of the discrete-time system~\eqref{eq_monotone_system-SGT}: there are $M>0$ and $a\in(0,1)$ such that for all $s \in\ell_\infty^+$ and $k\in\N_0$
\begin{equation*}
  \|\Gamma^k(s)\|_{\ell_\infty} \leq M a^k \|s\|_{\ell_\infty}.%
\end{equation*}
\item \emph{Uniform global asymptotic stability (UGAS)} of the discrete-time system~\eqref{eq_monotone_system-SGT}: there is $\beta\in\KL$ such that
\begin{equation*}
  \|\Gamma^k(s)\|_{\ell_\infty} \leq \beta(\|s\|_{\ell_\infty},k), \quad s \in\ell_\infty^+,\ k\in\N_0.
\end{equation*}
\item There is a \emph{point of strict decay}, i.e., there are $\lambda \in (0,1)$ and $s^0 \in \inner(\ell_{\infty}^+)$ such that
\begin{equation}\label{eq:Point-of-strict-decay}
  \Gamma(s^0) \leq \lambda s^0.%
\end{equation}
\item[(v)] It holds that%
\begin{equation}
\label{eq:Spectral-radius-SGC-checking-general}
\|\Gamma^n(\unit)\|_{\ell_\infty}<1 \mbox{\quad for some\ } n \in \N.%
\end{equation}
\end{enumerate}
\end{proposition}

\begin{proof}
(i) $\Leftrightarrow$ (ii). Follows from~\cite[Prop.~16 and its proof (in this part the subadditivity is used)]{MKG20}.%

(ii) $\Rightarrow$ (iii). Clear.%

(iii) $\Rightarrow$ (ii). From UGAS of system~\eqref{eq_monotone_system-SGT} and the homogeneity of $\Gamma$, we have%
\begin{equation*}
  \|\Gamma^k(s)\|_{\ell_{\infty}} = \|s\|_{\ell_{\infty}} \Big\|\Gamma^k\Big(\frac{s}{\|s\|_{\ell_\infty}}\Big)\Big\|_{\ell_\infty} \leq \|s\|_{\ell_{\infty}}\beta(1,k)% 
\end{equation*}
for all $0 \neq s \in \ell_{\infty}^+$ and $k \in \N_0$. As $\beta\in\KL$, there is $k_*\in\N$ such that $\beta(1,k_*)<1$. For each $k\in\N_0$, take $p,q \in\N_0$ satisfying $k = pk_* + q$ with $q<k_*$. Then%
\begin{align*}
  \|\Gamma^k(s)\|_{\ell_\infty} 
	&\leq \|\Gamma^q(s)\|_{\ell_\infty} \beta(1,k_*)^p \\
	&\leq \|s\|_{\ell_\infty} \beta(1,q) \beta(1,k_*)^p
	\leq \|s\|_{\ell_\infty} \beta(1,0) \beta(1,k_*)^p .
\end{align*}
Hence,~(ii) holds with $M = \beta(1,0)$ and $a = \beta(1,k_*)$.%

(iv) $\Rightarrow$ (ii). Pick any $s \in \ellinfp$. Then $s \leq \|s\|_{\ell_{\infty}} \unit \leq \frac{\|s\|_{\ell_{\infty}}}{s^0_{\inf}}s^0$, where
$s^0_{\inf} := \inf_{i\in\N} s_i^0$. As $s^0 \in \inner(\ellinfp) $, we have $s^0_{\inf}>0$. Since $\Gamma$ is monotone and homogeneous of degree one, we obtain%
\begin{equation*}
  \Gamma(s) \leq \Gamma\Big(\frac{\|s\|_{\ell_{\infty}}}{s^0_{\inf}}s^0\Big) = \frac{\|s\|_{\ell_{\infty}}}{s^0_{\inf}}\Gamma(s^0) \leq \frac{\|s\|_{\ell_{\infty}}}{s^0_{\inf}} \lambda s^0,%
\end{equation*}
and thus for all $k\in\N_0$ we have $\Gamma^k(s) \leq \frac{\|s\|_{\ell_{\infty}}}{s^0_{\inf}} \lambda^k s^0$. Due to monotonicity of the norm in $\ell_\infty^+$, we obtain~(ii):%
\begin{equation*}
  \|\Gamma^k(s)\|_{\ell_{\infty}} \leq \frac{\|s^0\|_{\ell_{\infty}}}{s^0_{\inf}} \lambda^k\|s\|_{\ell_{\infty}}.%
\end{equation*}

(ii) $\Rightarrow$ (iv). Pick any $\lambda \in (a,1)$, any $y \in \inner(X^+)$, and (motivated by \cite[Prop.~3.9]{GlM20}) consider the vector 
$z := \sum_{k=0}^\infty \frac{\Gamma^k(y)}{\lambda^{k+1}}$. The UGES property implies%
\begin{equation*}
  \sum_{k=0}^\infty \Big\|\frac{\Gamma^k(y)}{\lambda^{k+1}}\Big\|_{\ell_{\infty}} \leq \frac{M}{\lambda}\sum_{k=0}^\infty \frac{a^k}{\lambda^{k}}\|y\|_{\ell_{\infty}} <\infty.%
\end{equation*}
As $\ell_\infty$ is a Banach space, absolute convergence implies convergence, and thus $z$ is well-defined. As $\Gamma$ is monotone, we also have $z \geq \frac{1}{\lambda}y$, implying $z \in \inner(\ell_{\infty}^+)$ as $y \in \inner(X^+)$. Now decompose $\Gamma(z)$ as%
\begin{align*}
  \Gamma(z) &= \big(\Gamma - \lambda \id + \lambda \id\big)(z)\\
	&= \Gamma\Big( \sum_{k=0}^\infty \frac{\Gamma^k(y)}{\lambda^{k+1}}\Big) - \sum_{k=0}^\infty \frac{\Gamma^k(y)}{\lambda^{k}} + \lambda z.%
\end{align*}
Via subadditivity, continuity and homogeneity, we obtain%
\begin{align*}
  \Gamma(z) &\leq \sum_{k=0}^\infty \Gamma\Big(\frac{\Gamma^k(y)}{\lambda^{k+1}}\Big) - \sum_{k=0}^\infty \frac{\Gamma^k(y)}{\lambda^{k}} + \lambda z \\
&= \sum_{k=0}^\infty \frac{\Gamma^{k+1}(y)}{\lambda^{k+1}} - \sum_{k=0}^\infty \frac{\Gamma^k(y)}{\lambda^{k}} + \lambda z = -y + \lambda z \leq \lambda z,%
\end{align*}
which implies the claim.%

(i) $\Leftrightarrow$ (v). Using the subadditivity property of $\Gamma$, it follows that the limit on the right-hand side of the equality in~\eqref{eq_def_spectralradius} can be replaced by the infimum over all $n \in \N$:%
\begin{equation*}
  r(\Gamma) = \inf_{n \in \N} \sup_{s \in \ell_{\infty}^+,\ \|s\|_{\ell_\infty} = 1} \|\Gamma^n (s)\|_{\ell_{\infty}}^{1/n} =
	\inf_{n \in \N} \|\Gamma^n (\unit)\|_{\ell_{\infty}}^{1/n}.%
\end{equation*}
Hence, the condition $r(\Gamma) < 1$ is equivalent to \eqref{eq:Spectral-radius-SGC-checking-general}.
\end{proof}

As seen in Theorem~\ref{thm_mainres} below, Proposition~\ref{prop:eISS-criterion-homogeneous-systems} provides an equivalence between the small-gain condition for the verification of ISS of the infinite network~\eqref{eq_Sigma} and exponential stability of the discrete-time monotone system~\eqref{eq_monotone_system-SGT}. 
Moreover, condition~\eqref{eq:Point-of-strict-decay} yields a path of strict decay~\cite[Def.~II.10]{KMZ21} by which one constructs the overall Lyapunov function for the network $\Sigma$ as a supremum over weighted Lyapunov functions of subsystems $\Sigma_i$, where each weight is expressed as $1/s^0_i$ with $s^0$ as in~\eqref{eq:Point-of-strict-decay}. 

\subsection{Graph-theoretic aspects of the spectral radius condition} 
\label{sec:Computational aspects}

The condition \eqref{eq:Spectral-radius-SGC-checking-general} can be simplified, if $\Gamma$ has a special structure. If $\Gamma = \Gamma_\otimes$, the property \eqref{eq:Spectral-radius-SGC-checking-general} is equivalent to (see \cite[Lem.~12]{MKG20})
\begin{equation}\label{eq:Spectral-radius-SGC-checking}
  \sup_{j_1,\ldots,j_{n+1}} \gamma_{j_1j_2} \cdots \gamma_{j_nj_{n+1}} < 1 \mbox{\quad for some\ } n \in \N.%
\end{equation}
If $\Gamma = \Gamma_\oplus$, by induction, one can see that
\begin{equation*}
  \Gamma^n_\oplus(s) = \Big(\sum_{j_1,\ldots,j_{n-1}}\gamma_{i j_1}\gamma_{j_1j_2}\cdots\gamma_{j_{n-1}j} s_j \Big)_{i,j\in\N}.%
\end{equation*}
Hence, the property \eqref{eq:Spectral-radius-SGC-checking-general} is equivalent to
\begin{equation}\label{eq:Spectral-radius-SGC-checking-sum}
   \sup_{i\in\N}  \sum_{j_1,\ldots,j_n}\hspace{-3mm}\gamma_{i j_1}\gamma_{j_1j_2}\cdots\gamma_{j_{n-1}j_n} < 1 \mbox{\quad for some\ } n \in \N.%
\end{equation}
 
Let us provide a graph-theoretic description of \eqref{eq:Spectral-radius-SGC-checking} and~\eqref{eq:Spectral-radius-SGC-checking-sum}. For the infinite matrix $G := \big(\gamma_{ij}\big)_{i,j\in\N}$, we define a weighted directed graph graph ${\cal G}(G)$ with infinitely many nodes and the set of edges ${\cal E} = \{ (i,j) | i,j \in\N, \gamma_{ij} > 0\}$; thus, $(i,j) \in \cal E$ represents the edge from node $i$ to node $j$, and $G := \big(\gamma_{ij}\big)_{i,j\in\N}$ is the weight matrix with $\gamma_{ij}$ being the weight on edge $(i,j) \in \cal E$.
A path of length $n$ from $i$ to $j$ is a sequence of $n + 1$ distinct nodes starting with $i$ and ending with $j$.
Given ${\cal G}(G)$, condition~\eqref{eq:Spectral-radius-SGC-checking} requires that the product of weights over any path of length $n$ has to be uniformly less than one for some $n \in\N$.
For each node $i$ of ${\cal G}(G)$, condition~\eqref{eq:Spectral-radius-SGC-checking-sum} means that the sum over the product of weights over any path of length $n$ starting from node $i$ has to be uniformly less than one.
We stress that the two criteria~\eqref{eq:Spectral-radius-SGC-checking} and~\eqref{eq:Spectral-radius-SGC-checking-sum} are, in general, not comparable, as one is obtained from a max-formulation and the other from a sum-formulation, and thus their corresponding gains $\gamma_{ij}$ can be different. While, in general, it can be very hard (if not impossible) to verify such conditions for an infinite network, one can easily check conditions~\eqref{eq:Spectral-radius-SGC-checking} and~\eqref{eq:Spectral-radius-SGC-checking-sum} if the network has a special interconnection topology, e.g., is symmetric. We illustrate the latter point in Section~\ref{sec_ex} below, where we consider a spatially invariant network.

%%%%%%%%%%%%%%%%%%%%%%%%%%%%%%%

\subsection{The finite-dimensional case}

Although Proposition~\ref{prop:eISS-criterion-homogeneous-systems} was derived for monotone homogeneous subadditive operators acting in $\ell_\infty^+$, the same result also holds for operators acting in $\R^n_+$. 
Particularly, the results of Proposition~\ref{prop:eISS-criterion-homogeneous-systems} are novel even in finite-dimensional spaces and partially recover~\cite[Thm.~2]{Lu05} for $\Gamma = \Gamma_\otimes$ and~\cite[Lem.~1.1]{Rue10} for $\Gamma = \Gamma_\oplus$. Moreover, we have:%

\begin{proposition}
\label{prop:UGES-criterion-homogeneous-systems-finite-dim} 
Let an operator $\Gamma:\R^n_+\to\R^n_+$ be continuous, monotone, subadditive and homogeneous of degree one. Then $r(\Gamma)<1$ if and only if there is a \emph{path of strict decay with respect to $\Gamma$}, i.e., if there are maps $\sigma =(\sigma_i)_{i=1}^n \in\Kinf^n$ and $\rho \in\Kinf$ such that%
\begin{eqnarray}\label{eq:Path-of-strict-decay-inequality}
  \Gamma(\sigma(r)) \leq (\id + \rho)^{-1}(\sigma(r)) \mbox{\quad for all\ } r \geq 0.%
\end{eqnarray}
\end{proposition} 

\begin{proof}
$\Rightarrow$. The condition $r(\Gamma)<1$ implies the existence of $\lambda \in (0,1)$ and $s^0 \in \inner(\ell_{\infty}^+)$ such that $\Gamma(s^0) \leq \lambda s^0$, and thus $\Gamma(\sigma(r)) \leq \lambda \sigma(r)$ for all $r\geq 0$, where $\sigma(r):=rs^0$. Thus, $\sigma$ is a linear path of strict decay.%

$\Leftarrow$. The existence of a path of strict decay with respect to $\Gamma$ implies UGAS of the corresponding discrete-time system due to \cite[Thm.~4]{RuS14} with $A(s):=g(s,0)$.
Item (iii) of Proposition~\ref{prop:eISS-criterion-homogeneous-systems} completes the proof.
\end{proof}

%%%%%%%%%%%%%%%%%%%%%%%%%%%%%%%

\section{Small-gain theorem} \label{sec:Small-gain theorem}

Now we are able to present our small-gain theorem which provides a Lipschitz continuous ISS Lyapunov function for the infinite interconnection $\Sigma$.%

\begin{theorem}\label{thm_mainres}
Consider the infinite network $\Sigma$ as in \eqref{eq_Sigma}, and let the following assumptions hold:%
\begin{enumerate}[(i)]
%\item The network $\Sigma$ is well-posed and satisfies the BIC property.%
\item $\Sigma$ is well-posed and satisfies the BIC property.%
\item Assumption \ref{ass_subsystem_iss} is satisfied. Moreover, there exist $\psi_1,\psi_2 \in \KC_{\infty}$ such that%
\begin{equation}\label{eq_coerc_unif_bounds}
  \psi_1 \leq \psi_{i1} \mbox{\quad and \quad} \psi_{i2} \leq \psi_2 \mbox{\quad for all\ } i \in \N.%
\end{equation}
\item Assumptions \ref{ass_gainop_wd} and \ref{ass_gammaiu_max} are satisfied.%
\item The spectral radius condition $r(\Gamma) < 1$ holds.%
\item For each $R>0$, there is $L(R)>0$ such that%
\begin{equation}\label{eq_lyap_lipschitz}
  |V_i(x_i) - V_i(y_i)| \leq L(R)|x_i - y_i|%
\end{equation}
for all $i\in\N$ and $x_i,y_i \in B_R(0) \subset \R^{n_i}$.
\item There exists $\tilde{\alpha} \in \PC$ such that $\alpha_i \geq \tilde{\alpha}$ for all $i\in\N$.%
\end{enumerate}
Then there exists $s^0 \in \inner(\ellinfp)$ such that the following function is an ISS Lyapunov function for $\Sigma$:%
\begin{equation}\label{eq_def_V}
  V(x) := \sup_{i\in\N} \frac{1}{s^0_i} V_i(x_i) \mbox{\quad for all\ } x = (x_i)_{i\in\N} \in X.%
\end{equation}
In particular, $\Sigma$ is ISS. Moreover, $V$ satisfies a Lipschitz condition on every bounded subset of $X$.%
\end{theorem}

Large parts of the proof of Theorem \ref{thm_mainres} are analogous to corresponding parts of the proof of~\cite[Thm.~III.1]{KMZ21}, and hence will be omitted.
\begin{proof}
By Proposition~\ref{prop:eISS-criterion-homogeneous-systems} and since $r(\Gamma)<1$, we find $s^0 \in \inner(\ellinfp)$ and $\lambda\in(0,1)$ with $\Gamma(s^0) \leq \lambda s^0$. Since $s^0$ is an interior point of the cone $\ell_{\infty}^+$, there are $s^0_{\min},s^0_{\max}>0$ such that the components $s^0_i$ of $s^0$ satisfy $s^0_{\min} \leq s^0_i \leq s^0_{\max}$ for all $i \in \N$, which is used throughout the proof. From \eqref{eq_subsystem_iss_coerc} and \eqref{eq_coerc_unif_bounds}, it easily follows that $V$ satisfies the estimates%
\begin{equation}\label{eq_coercivity_est}
  \frac{1}{s^0_{\max}} \psi_1(\|x\|_X) \leq V(x) \leq \frac{1}{s^0_{\min}} \psi_2(\|x\|_X),\quad \forall x \in X.
\end{equation}
Also, for any $R>0$ and $x,y \in B_{R}(0)$, we have $|V(x) - V(y)| \leq \frac{1}{s^0_{\min}}L(R)\|x-y\|_X$. Now, we fix a number $\mu \in (1,\lambda^{-1})$ and introduce, for every $0 \neq x \in X$, the set%
\begin{equation*}
  I(x) := \Bigl\{ i \in \N : V(x) \leq \frac{\mu}{s^0_i}V_i(x_i) \Bigr\}.%
\end{equation*}
Then, with the same arguments as used in \cite[Proof of Thm.~III.1]{KMZ21}, we can show that for every $0 \neq x \in X$ there is $\delta>0$ such that%
\begin{equation}\label{eq_ixrep}
  V(y) = \sup_{i \in I(x)} \frac{1}{s^0_i}V_i(y_i) \mbox{\quad for all\ } y \in B_{\delta}(x) .
\end{equation}
Define $\gamma\in\KC$ by $\gamma(r) := \frac{1}{s^0_{\min}\lambda} \gamma_{\max}^u(r)$ for all $r \geq 0$, and fix $0 \neq x \in X$ and $u \in \UC$ satisfying $V(x) > \gamma(\|u\|_{\UC})$.
We first show that there is $T>0$ such that%
\begin{equation}\label{eq_iss_subsystem_est}
  \nabla V_i(\phi_i(t)) f_i(\phi_i(t),\bar{\phi}_i(t),u_i(t)) \leq -\tilde{\alpha}(V_i(\phi_i(t)))%
\end{equation}
for all $i \in I(x)$ and $t \in [0,T]$, where $\phi(t) := \phi(t,x,u)$, $\phi_j(t)$ is the $j$-th component of $\phi(t)$ (for every $j\in\N$) and $\bar{\phi}_i(t) = (\phi_j(t))_{j\in I_i}$. To prove \eqref{eq_iss_subsystem_est}, note that the continuity of $V$ and $\phi(\cdot)$ implies%
\begin{equation}\label{eq_lyap_hypothesis}
  V(\phi(t)) > \frac{1}{s^0_{\min}\lambda} \gamma_{\max}^u(\|u\|_{\UC}) \geq \frac{1}{s^0_{\min}\lambda}\gamma_{iu}(|u_i(t)|)%
\end{equation}
for all $i \in \N$, $t \in [0,T]$, where $T>0$ is chosen sufficiently small.
If $T$ is chosen further small enough, then%
\begin{equation}\label{eq_subclaim}
  V(\phi(t)) < \frac{\lambda^{-1}}{s^0_i}V_i(\phi_i(t)) \mbox{\quad for all\ } i \in I(x),\ t \in [0,T],%
\end{equation}
which is proved as in \cite{KMZ21} (Step 4 of the proof of Thm.~III.1 therein). From \eqref{eq_subclaim}, it then follows that%
\begin{align*}
  V_i(\phi_i(t)) &> V(\phi(t)) \lambda s^0_i \geq V(\phi(t))\pi_i \circ \Gamma_{\mu}(s^0) \\ 
	               &= V(\phi(t)) \mu_i( (\gamma_{ij} s^0_j)_{j\in\N} ) = \mu_i( (\gamma_{ij} V(\phi(t)) s^0_j)_{j\in\N} ) \\
	               &\geq \mu_i( (\gamma_{ij} V_j(\phi_j(t)))_{j\in\N} ),%
\end{align*}
where we use the inequality $\Gamma(s^0) \leq \lambda s^0$ and the assumption that $\mu_i$ is homogeneous and monotone. At the same time, \eqref{eq_lyap_hypothesis} together with \eqref{eq_subclaim} implies $V_i(\phi_i(t)) > \gamma_{iu}(|u_i(t)|)$. Putting both estimates together, we obtain%
\begin{equation*}
  V_i(\phi_i(t)) > \max\Bigl\{ \mu_i \Big(\big(\gamma_{ij} V_j (\phi_j(t))\big)_{j\in\N}\Big) ,\gamma_{iu}(|u_i(t)|)\Bigr\}.%
\end{equation*}
By \eqref{eq_subsystem_orbitalder_est} together with Assumption (vi), this implies \eqref{eq_iss_subsystem_est}. Hence, the proof of \eqref{eq_iss_subsystem_est} is complete. We now consider the scalar differential equation%
\begin{equation}\label{eq:comparison-2}
  \dot{v}(t) = -\tilde{\alpha}(v(t)), \quad v(0) = v_0 \in \R_+.%
\end{equation}
By \cite[Lem.~4]{MKG20}, we can assume that $\tilde{\alpha}$ satisfies a global Lipschitz condition, so that the Cauchy problem \eqref{eq:comparison-2} has a unique and globally defined solution $\VC(t,v_0)$, $t \geq 0$. It follows from a standard comparison lemma, e.g.,~\cite[Lem.~4.4]{LSW96}, that for every $i \in I(x)$ and all $t$ sufficiently small, the inequality $V_i(\phi_i(t)) \leq \VC(t,V_i(x_i))$ is satisfied. Hence, for all $t>0$ sufficiently small, one has
\begin{align*}
   \frac{1}{t}(V(\phi(t)) &- V(x)) = \frac{1}{t}\Bigl(\sup_{i\in I(x)} \frac{1}{s^0_i}V_i(\phi_i(t)) \!-\! \!\sup_{i \in I(x)} \!\frac{1}{s^0_i}V_i(x_i)\Bigr) \\
	&\leq \frac{1}{t}\sup_{i\in I(x)}\Bigl(\frac{1}{s^0_i}V_i(\phi_i(t)) - \frac{1}{s^0_i}V_i(x_i)\Bigr) \\
	&\leq	\sup_{i \in I(x)} \frac{1}{t} \frac{1}{s^0_i}\left(\VC(t,V_i(x_i)) - V_i(x_i)\right) \\
	&= - \inf_{i \in I(x)} \frac{1}{t} \frac{1}{s^0_i} \int_0^t \tilde{\alpha}(\VC(s,V_i(x_i)))\, \rmd s \\
	&\leq -\frac{1}{s^0_{\max}} \frac{1}{t}\int_0^t \min_{\zeta \in [\mu^{-1}s^0_{\min},s^0_{\max}]} \tilde{\alpha}(\VC(s,\zeta V(x)))\, \rmd s.%
\end{align*}
As $s \mapsto \min_{\zeta \in [\mu^{-1} s^0_{\min},s^0_{\max}]} \tilde{\alpha}(\VC(s,\zeta V(x)))$ is a continuous function, letting $t \rightarrow 0^+$ leads to%
\begin{equation*}
  \rmD^+ V_u(x) \leq -\frac{1}{s^0_{\max}}\min_{\zeta \in [\mu^{-1}s^0_{\min},s^0_{\max}]} \tilde{\alpha}(\zeta V(x)).%
\end{equation*}
It is easy to see that $\alpha(r) := \frac{1}{s^0_{\max}} \min_{\zeta \in [\mu^{-1} s^0_{\min},s^0_{\max}]} \tilde{\alpha}(\zeta r)$ is a positive definite function. Hence, we have
\begin{equation}
\label{eq_limp}
  V(x) > \gamma(\|u\|_{\UC}) \quad \Rightarrow \quad \rmD^+ V_u(x) \leq -\alpha(V(x))%
\end{equation}
which completes the proof.%
\end{proof}

Theorem~\ref{thm_mainres} shows ISS of the network $\Sigma$ under the small-gain condition $r(T)<1$, which by Proposition~\ref{prop:eISS-criterion-homogeneous-systems} is equivalent to the existence of a point of strict decay, and thus to the existence of a linear path of strict decay with a linear decay rate. Note that in \cite{KMZ21} ISS of the network has been shown under the assumption of the existence of a nonlinear path of strict decay with a nonlinear decay rate. Whether these conditions are equivalent for homogeneous of degree one and subadditive gain operators is an open question. For finite networks, Proposition~\ref{prop:UGES-criterion-homogeneous-systems-finite-dim} shows the equivalence of these two concepts.

%%%%%%%%%%%%%%%%%%%%%%%%%%%%%%%%%%%%%%%%%%%%%%
\subsection{Dissipative formulation}\label{sec:Dissipative formulation}

So far, we have only considered ISS Lyapunov functions in an implication form.
A dissipative formulation of ISS Lyapunov functions is also quite useful in various applications. A central question in the context of ISS is whether these two Lyapunov formulations are equivalent. It is not hard to see that the implication from an ISS Lyapunov function in dissipative form to one in implication form holds for both finite- and infinite-dimensional systems. For infinite-dimensional systems, however, it is an open question whether, in general, the converse holds, see~\cite{MiP20}. Nevertheless, here we provide additional, but mild conditions under which the converse holds, at least within the setup of this paper. First recall the notion of an ISS Lyapunov function in a dissipative form.

\begin{definition}\label{def:dissipative-form}
A function $V:X \rightarrow \R_+$ is called an \emph{ISS Lyapunov function (in a dissipative form)} for $\Sigma$ if it satisfies items (i) and (ii) of Definition~\ref{def:implication-form}, and additionally, there are $\alpha \in \KC_{\infty}$ and $\rho \in \KC$ such that for all $x\in X$ and $u \in \UC$ we have
\begin{equation}\label{eq_orbitalder_diss_est}
  \rmD^+ V_u(x) \leq -\alpha(V(x)) + \rho(\|u\|_{\UC}).%
\end{equation}
\end{definition}

A notion of \emph{nonlinear scaling} is required to present the result of this section.
\begin{definition}\label{def:nonlinear scaling} 
A \emph{nonlinear scaling} is a function $\xi\in\Kinf$ which is continuously differentiable on $(0,\infty)$ and satisfies $\xi'(s)>0$ for all $s>0$ such that the limit $\xi'(0) := \lim_{s\to 0}\xi'(s)$ exists and is finite.%
\end{definition}

Now we present the result of this section.%

\begin{proposition}
\label{prop:Nonlinear-scaling}
Additionally to the assumptions of Theorem \ref{thm_mainres}, let the following hold:%
\begin{enumerate}[(i)]
\item The function $f:X \tm U \rightarrow X$ is continuous and for every $R>0$ it holds that%
\begin{equation}\label{eq_extra_ass_diss2}
  C(R) := \sup_{\max\{\|x\|_X,\|u\|_U\} \leq R} \|f(x,u)\|_X < \infty.%
\end{equation}
\item The ISS Lyapunov functions $V_i$, $i\in\N$, are continuously differentiable at the origin and the constants $L(R)$ in \eqref{eq_lyap_lipschitz} can be chosen as continuous functions of $R$ satisfying $L(R) \rightarrow 0$ as $R \rightarrow 0$.%
\end{enumerate}
Then there exists a nonlinear scaling $\xi\in\Kinf$, such that for $V$ defined in~\eqref{eq_def_V}, the function $W := \xi \circ V$ is an ISS Lyapunov function for $\Sigma$ in dissipative form.
\end{proposition}

\begin{proof}
Let $0 \neq x \in X$, $u \in \UC$ and write $\phi(t) := \phi(t,x,u)$. We start with an estimate on $\rmD^+ V_u(x)$, using the mean value theorem.
For any $t>0$ sufficiently small, we have%
\begin{align*}
&  \frac{1}{t}\left(V(\phi(t)) - V(x)\right) \leq \sup_{i\in\N} \frac{1}{t} \frac{1}{s^0_i}(V_i(\phi_i(t)) - V_i(x_i)) \\
&\quad = \sup_{i\in\N} \frac{1}{s^0_i} \nabla V_i(\phi_i(t_i)) f_i(\phi_i(t_i),\bar{\phi}_i(t_i),u_i(t)) \\
&\quad \leq \frac{1}{s^0_{\min}} \sup_{i\in\N} |\nabla V_i(\phi_i(t_i))| |f_i(\phi_i(t_i),\bar{\phi}_i(t_i),u_i(t))|,
\end{align*}
with $t_i \in (0,t)$. Since $|\phi_i(t_i)| \leq \|\phi(t_i)\|_X$, with $R_{x,u}(t) := \max_{s\in[0,t]} \|\phi(s)\|_X$, we obtain%
\begin{equation*}
  \frac{1}{t}\left(V(\phi(t)) - V(x)\right) \leq \frac{1}{s^0_{\min}} L(R_{x,u}(t)) \sup_{s\in[0,t]} \|f(\phi(s),u(s))\|_X,%
\end{equation*}
where we use that the norm of the gradient can be bounded by the Lipschitz constant. By continuity of $f$, $\phi(\cdot)$ and $u(\cdot)$ (in a neighborhood of $0$) and Assumption (ii), we obtain%
\begin{equation*}
  \rmD^+ V_u(x) \leq \frac{1}{s^0_{\min}} L(\|x\|_X) \|f(x,u(0))\|_X.%
\end{equation*}
Let the functions $\gamma$ and $\alpha$ be constructed as in the proof of Theorem \ref{thm_mainres} (see \eqref{eq_limp}). As the decay rate $\alpha$ in \eqref{eq_limp} may be a positive definite function only, following \cite[Rem.~4.1]{LSW96}, one can find a nonlinear scaling $\xi\in\Kinf$ of the form $\xi(r) = \int_0^r\frac{1}{\sigma(s)}\, \rmd s$, $r\geq 0$, for some $\sigma\in\LL$, and $\bar{\alpha}\in\Kinf$ such that $W = \xi\circ V$ is an ISS Lyapunov function satisfying
\begin{equation}\label{ImplicationIneq_cont-aux3}
  W(x) > \xi \circ \gamma(\|u\|_\UC) \quad \Rightarrow \quad \rmD^+ W_u(x) \leq -\bar{\alpha}(W(x)).
\end{equation}
It holds that%
\begin{align*}
  \rmD^+ W_u(x) & = \frac{1}{\sigma(V(x))} \rmD^+ V_u(x) \\
  & \leq \frac{1}{\sigma(V(x))} \frac{1}{s^0_{\min}} L(\|x\|_X) \|f(x,u(0))\|_X .
\end{align*}
Define $\tilde{\rho}(r) := \sup\Bigl\{ \rmD^+ W_u(x) + \bar\alpha(\xi\circ\gamma(\|u\|_{\UC})) : \|u\|_{\UC} \leq r,\ W(x) \leq \xi\circ\gamma(\|u\|_{\UC}) \Bigr\}$ 
and $\rho(r) := \max\{0,\tilde{\rho}(r)\}$. Let us also put $
  \|(x,u)\| := \max\{\|x\|_X,\|u\|_{\UC}\} \mbox{\ and\ } R(r) := \max\{r,\psi_1^{-1}(s^0_{\max}\gamma(r))\}$. 
As $W(x) \leq \xi\circ\gamma(\|u\|_{\UC})$ is equivalent to $V(x) \leq \gamma(\|u\|_{\UC})$, with $\|x\|_X \leq \psi_1^{-1}(s^0_{\max}V(x))\leq \psi_1^{-1}(s^0_{\max}\gamma(\|u\|_{\UC}))$ (see \eqref{eq_coercivity_est}) we obtain 
\begin{align*}
  \rho(r)& \leq  \frac{1}{\sigma(\gamma(r))} \frac{1}{s^0_{\min}} L(\psi_1^{-1}(s^0_{\max} \gamma(r))) \!\!\!\!
  \sup_{\|(x,u)\| \leq R(r)} \|f(x,u)\|_X \\
  &\quad + \bar\alpha(\xi\circ\gamma(r) )
	\\
	&= \frac{1}{\sigma(\gamma(r))} \frac{1}{s^0_{\min}} L(\psi_1^{-1}(s^0_{\max}\gamma(r))) C(R(r))
	 + \bar\alpha(\xi\circ\gamma(r)).%
\end{align*}
In particular, $\rho$ is nondecreasing, $\rho(0) = 0$ and $\lim_{r \rightarrow 0^+}\rho(r) = 0$ by Assumption (ii). 
One may assume that $\rho \in \KC$ (otherwise, majorize $\rho$ by a $\KC$-function, which is always possible). To show \eqref{eq_orbitalder_diss_est}, we distinguish two cases:\\
- Assume that $W(x) > \xi\circ \gamma(\|u\|_{\UC})$.
Then~\eqref{eq_limp} yields%
\begin{equation*}
   \rmD^+ W_u(x) \leq -\bar\alpha(W(x)) \leq -\bar\alpha(W(x)) + \rho(\|u\|_{\UC}).%
\end{equation*}
- Assume that $W(x) \leq \xi\circ \gamma(\|u\|_{\UC})$. Then%
\begin{align*}
  \rmD^+ W_u(x) &= \rmD^+ W_u(x) + \bar\alpha(\xi\circ\gamma(\|u\|_{\UC})) - \bar\alpha(\xi\circ\gamma(\|u\|_{\UC})) \\
	              &\leq \tilde{\rho}(\|u\|_{\UC}) - \bar\alpha(\xi\circ\gamma(\|u\|_{\UC})) \\
	              &\leq -\bar\alpha(W(x)) + \rho(\|u\|_{\UC}).
\end{align*}
This concludes the proof.
\end{proof}

%%%%%%%%%%%%%%%%%%%%%%%%%%%%%%%%%%%%%%%%%%%%%%

\vspace{-0.3cm}
\section{Example}\label{sec_ex}

In this section, we illustrate the effectiveness of conditions~\eqref{eq:Spectral-radius-SGC-checking} and~\eqref{eq:Spectral-radius-SGC-checking-sum} to verify the small-gain condition.
In particular, we allow the gain functions in part to be larger than the identity.%

Consider scalar linear subsystems of the form%
\small
\begin{equation*}
  \Sigma_i: \! \left\{\begin{array}{l}
\!\!\dot{x}_i = - b_{ii}x_i + b_{i (i-1)}x_{i-1}+ g (b_{i(i+1)}x_{i+1},b_{i(i+2)}x_{i+2}) \\
\qquad\qquad\qquad\qquad\qquad\qquad\,\textrm{for }  i = 2k+1,\ k\in\N ,\\
\!\!\dot{x}_i = - b_{ii}x_i  + g(b_{i(i+1)}x_{i+1},b_{i(i+2)}x_{i+2}) \\
 \qquad\qquad\qquad\qquad\qquad\qquad\,\textrm{ for } i = 2k,\ k\in\N,
\end{array}\right.%
\end{equation*}
\normalsize{
$\!\!$where either $g(s) =: g_\oplus(s) = s_1+s_2$ or $g(s) =: g_\otimes(s) = \max\{s_1,s_2\}$ for $s = (s_1,s_2) \in\R^2$, $x_i \in\R$, $b_{ii} > 0$ and $b_{i(i-1)},b_{i(i+1)},b_{i(i+2)} \in \R$ with $b_{i0} = 0$.}

%The structure of the network is depicted in Figure~\ref{topology}.%

%\begin{figure}
%	\centering  
%	\includegraphics[totalheight=1.2cm]{topology-2}
%	\caption{The interconnection topology of the network in Section \ref{sec_ex}.}
%	\label{topology}
%\end{figure}

To guarantee that the network $\Sigma$ of subsystems $\Sigma_i$ is well-posed and satisfies the BIC property, it suffices to assume that all the coefficients $b_{ij}$ are uniformly upper-bounded by some $\overline{b}>0$ over $i$ (the simple proof will be omitted).

     \begin{figure}[h]
     \vspace{-0.3cm}
\centering
%\hspace{-2cm}
    \includegraphics[trim={5.2cm 10.2cm 2.0cm 0cm},clip,width=\linewidth]{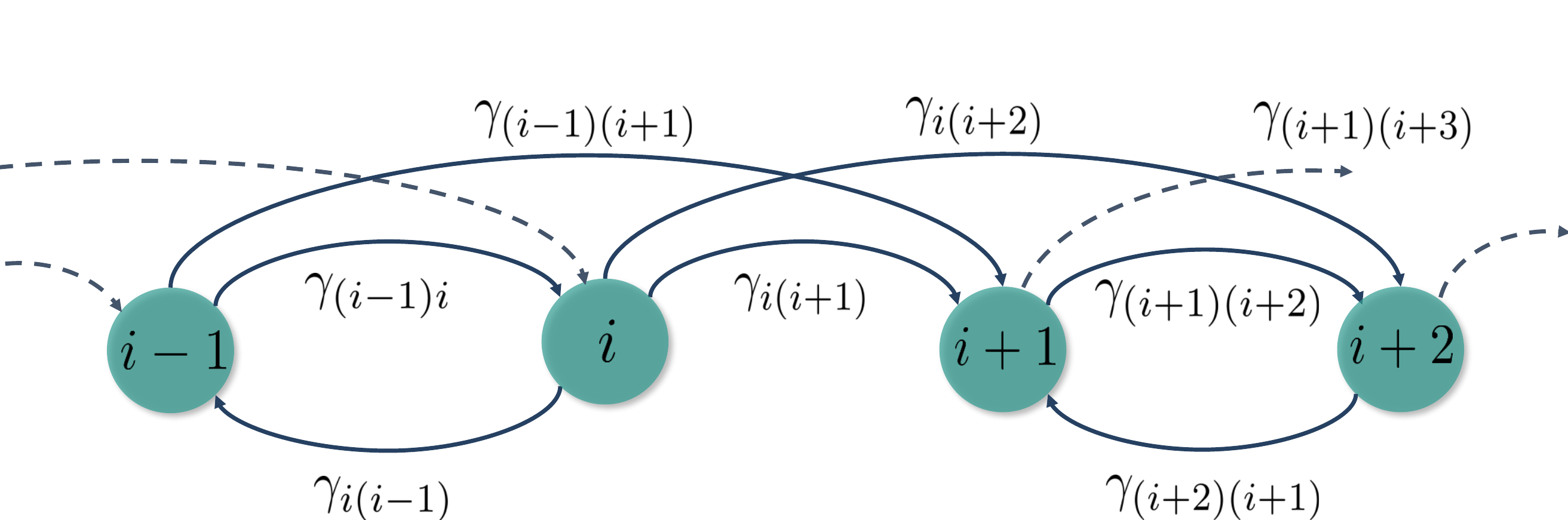}
\vspace{-0.6cm}
\caption{The graph ${\cal G}(G)$ associated with $G$. Dashed arrows imply the continuation of the graph with the same structure.}
\label{fig:topology}
\end{figure}
Take $V_i(x_i) := \frac{1}{2}x_i^2$. 
The following estimates are obtained by using Young's inequality and homogeneity of $g$:%
\begin{align}\label{eq_linex_diss_est}
\begin{split}
&  \nabla V_i(x_i) f_i(x_i,\bar{x}_i) \!\leq\! - \alpha_i V_i(x_i) \!+\! \frac{b_{i(i-1)}^2}{2\ep_i}V_{i-1}(x_{i-1}) \\
	&+g\big( \frac{b_{i(i+1)}^2}{2\delta_i}V_{i+1}(x_{i+1}) , \frac{b_{i(i+2)}^2}{2\delta^\prime_i}V_{i+2}(x_{i+2})\big) ,
\end{split}
\end{align}
with $\alpha_i := 2(b_{ii}-\ep_i-\delta_i-\delta_i^\prime)$, and appropriately small $\delta_i,\delta^\prime_i,\ep_i>0$. If $i$ is even, then $b_{i(i-1)} = 0$ and neither the second term nor $\ep_i$ on the right-hand side of~\eqref{eq_linex_diss_est} is present. Additionally, if $g = g_\otimes$, then $\delta_i^\prime = \delta_i$ in $g$ and $\alpha_i := 2(b_{ii}-\ep_i-\delta_i)$. From this estimate, we need to derive appropriate gains for an implication formulation of the ISS property.
Let
 $
\tilde{\gamma}_{i(i-1)} := b_{i(i-1)}^2/2\ep_i$, $\tilde{\gamma}_{i(i+1)} := b_{i(i+1)}^2/2\delta_i$, $\tilde{\gamma}_{i(i+2)} := b_{i(i+2)}^2/2\delta^\prime_i
$
and assume that
%\begin{align}
%&V_i(x_i) \geq \frac{1}{a_i} \tilde\mu ( (\tilde{\gamma}_{ik} V_{k}(x_{k}))_{k=-1,1,2}) \label{eq_linex_ass}
%\end{align}
\begin{align}
&V_i(x_i) \geq \frac{1}{a_i} \tilde\mu \big( (\tilde{\gamma}_{i(i+k)} V_{i+k}(x_{i+k}))_{k=-1,1,2}\big)  \label{eq_linex_ass} ,
\end{align}
%\begin{align}
%&V_i(x_i) \geq \nonumber\\
%&\frac{1}{a_i} \tilde\mu ( \tilde{\gamma}_{i(i-1)} V_{i-1}(x_{i-1}), \tilde{\gamma}_{i(i+1)} V_{i+1}(x_{i+1}),\tilde{\gamma}_{i(i+2)} V_{i+2}(x_{i+2}) ) \label{eq_linex_ass}
%\end{align}
where in line with the function $g$ we have either $\tilde\mu(s) =  2 \max_{i=1,2,3}s_i$ or $\tilde\mu(s) = \sum_{i=1}^{3} s_i$, and $a_i > 0$ is specified later.
From~\eqref{eq_linex_diss_est}, it follows that
\begin{align*}
 & \nabla V_i(x_i) f_i(x_i,\bar{x}_i) \leq  -\alpha_iV_i(x_i) \\
 &+ \tilde\mu( \tilde{\gamma}_{i(i-1)} V_{i-1}(x_{i-1}), \!\tilde{\gamma}_{i(i+1)} V_{i+1}(x_{i+1}),\!\tilde{\gamma}_{i(i+2)} V_{i+2}(x_{i+2})).
\end{align*}
Hence,~\eqref{eq_linex_ass} implies 
$\nabla V_i(x_i) f_i(x_i,\bar{x}_i) \leq -(\alpha_i - a_i)V_i(x_i)$.
Let us pick $a_i := \alpha_i/2$.
In that way, for $g = g_\oplus$, non-zero gains $\gamma_{ij}$ in $\Gamma_\oplus$ are computed as
\begin{align*}
\gamma_{i(i-1)} = \frac{b_{i(i-1)}^2}{2\ep_i w_i}, \,\,
\gamma_{i(i+1)} = \frac{b_{i(i+1)}^2}{2\delta_i w_i},\,\,
\gamma_{i(i+2)} = \frac{b_{i(i+2)}^2}{2\delta_i^\prime w_i} ,
\end{align*}
where $w_i := b_{ii}-\ep_i-\delta_i-\delta_i^\prime$.
On the other hand, if $g = g_\otimes$, non-zero gains $\gamma_{ij}$ in $\Gamma_\otimes$ are obtained as
\begin{align*}
&\gamma_{i(i-1)} = \frac{b_{i(i-1)}^2}{\ep_i q_i } , \,
\gamma_{i(i+1)} = \frac{b_{i(i+1)}^2}{\delta_i q_i} ,\,
\gamma_{i(i+2)} = \frac{b_{i(i+2)}^2}{\delta_i q_i} ,
\end{align*}
with $q_i = b_{ii}-\ep_i-\delta_i$.
Obviously, Assumptions (i), (ii), and (v) of Theorem \ref{thm_mainres} are satisfied. To verify Assumption (vi), it is necessary to impose a uniform lower bound on the number $b_{ii}$, i.e., $b_{ii} \geq \underline{b}$ for all $i \in \N$.
Furthermore, we require uniform bounds on the numbers $\delta_i$, $\delta_i^\prime$, $\ep_i$ so that 
$\underline{\delta} \leq \delta_i \leq \overline{\delta}$, $\underline{\delta}^\prime \leq \delta_i^\prime \leq \overline{\delta}^\prime$, $\underline{\ep} \leq \ep_i \leq \overline{\ep}$, and $\underline{b} - \overline{\ep} - \overline{\delta}- \overline{\delta}^\prime > 0$.
Then for both sum and max formulations, we can estimate 
$
\gamma_{ij} \leq 2\ol{b}^2/\big(\underline{\ep}(\ul b - \ol{\ep} - \ol{\delta} - \ol{\delta}^\prime\big)
$ for all $i,j \in \N$.
Hence, Assumptions (iii) and (vi) of Theorem \ref{thm_mainres} are satisfied.

It only remains to verify the small-gain condition. 
To do this, we respectively use the criteria~\eqref{eq:Spectral-radius-SGC-checking} and~\eqref{eq:Spectral-radius-SGC-checking-sum} for $g = g_\otimes$ and $g = g_\oplus$.
Following arguments of Section~\ref{sec:Computational aspects}, we use define the weighted directed graph ${\cal G}(G)$ which is (partially) depicted in Fig.~\ref{fig:topology}.
As can be seen, the graph is symmetrical due to the spatial invariance of the overall network.
Take $n = 2$ in~\eqref{eq:Spectral-radius-SGC-checking}.
From the structure of the network and the choice of $n$, condition~\eqref{eq:Spectral-radius-SGC-checking} reduces to the following set of conditions:
%\begin{align}\label{eq:sgc-odd}
%\begin{array}{rl}
 %\gamma_{(i-1)(i+1)} \gamma_{(i+1)(i+3)} & < 1,\\
  %\gamma_{(i-1)(i+1)} \gamma_{(i+1)(i+2)} & < 1, \\
  %\gamma_{(i-1)i} \gamma_{i(i-1)} & < 1,\\
    %\gamma_{(i-1)i} \gamma_{i(i+1)} & < 1,\\
  %\gamma_{(i-1)i} \gamma_{i(i+2)} & < 1,
%\end{array}
%\end{align}
\begin{align}\label{eq:sgc-odd}
\begin{array}{rl}
 \gamma_{(i-1)(i+1)} \gamma_{(i+1)(i+k)} & < 1, \quad k\in\{2,3\},\\
  %\gamma_{(i-1)(i+1)} \gamma_{(i+1)(i+2)} & < 1, \\
  \gamma_{(i-1)i} \gamma_{i+k} & < 1, \quad k\in\{-1,1,2\},
    %\gamma_{(i-1)i} \gamma_{i(i+1)} & < 1,\\
  %\gamma_{(i-1)i} \gamma_{i(i+2)} & < 1,
\end{array}
\end{align}
for all $i \in\N$.
On the other hand, by taking node $i-1$, condition~\eqref{eq:Spectral-radius-SGC-checking-sum} yields
\begin{align}
& \!\!\!\!\!\gamma_{(i-1)(i+1)} \gamma_{(i+1)(i+3)} +
  \gamma_{(i-1)(i+1)} \gamma_{(i+1)(i+2)} \nonumber\\
&  \!\!\!+  \gamma_{(i-1)i} \gamma_{i(i-1)} +
    \gamma_{(i-1)i} \gamma_{i(i+1)} +
  \gamma_{(i-1)i} \gamma_{i(i+2)}  < 1 . \label{eq:sgc-odd-sum}
\end{align}
Clearly, any other nodes give the same condition as that in~\eqref{eq:sgc-odd-sum}.
From either~\eqref{eq:sgc-odd} or~\eqref{eq:sgc-odd-sum}, for sufficiently small $\gamma_{(i-1)i}$ and $\gamma_{(i-1)(i+1)}$, conditions~\eqref{eq:sgc-odd} and~\eqref{eq:sgc-odd-sum} are satisfied.
%One such choice is $\gamma_{(i-1)i} := 0.01$, $\gamma_{i(i-1)} := 2$, $\gamma_{i(i+1)} := 4$.
Given $b_{ii}, \ep_i, \delta_i,\delta_i^\prime$, as $\gamma_{ij}$ directly depend on $b_{ij}$, to satisfy any of~\eqref{eq:sgc-odd} and~\eqref{eq:sgc-odd-sum}, $b_{ij}$ should be small enough.

% $ \gamma_{(i-1)(i+1)} =  \gamma_{(i+1)(i+3)}$ that due to the spatial invariance of the network

%\textcolor{red}{[The way the example is written is not acceptable. Long formulas must be written in separated lines.]}

\section{Conclusions}

We developed a Lyapunov-based small-gain theorem ensuring ISS for infinite networks. 
Assuming linear internal gain functions, we particularly considered monotone, subadditive and homogeneous gain operators. 
Our small-gain condition is expressed in terms of the spectral radius of the gain operator. 
Several equivalences of the spectral radius condition were provided. 
In particular, we showed that the spectral radius condition is equivalent to UGES of the discrete-time system induced by the gain operator. 
Two special cases of the gain operator are linear and max-linear operators, for which we provided explicit formulas for the spectral radius. 
Computational aspects of the latter conditions were discussed via an illustrative example.

\bibliographystyle{IEEEtran}
\bibliography{references}

\end{document}